\documentclass[a4paper,UKenglish]{lipics}
\usepackage{microtype}%if unwanted, comment out or use option "draft"

\title{Teaching Logics through Their Philosophical Commitments: “Logical Worldviews”}
\titlerunning{Teaching Logics through Their Philosophical Commitments} %optional, in case that the title is too long; the running title should fit into the top page column

\author{Creighton Rosental}
\affil{Associate Professor of Philosophy\\
  Mercer University\\
  \texttt{rosental\_c@mercer.edu}}

\authorrunning{C. Rosental} %mandatory. First: Use abbreviated first/middle names. Second (only in severe cases): Use first author plus 'et. al.'

\Copyright{Creighton Rosental}%mandatory, please use full first names. LIPIcs license is "CC-BY";  http://creativecommons.org/licenses/by/3.0/

% mandatory: Please provide 1-5 keywords
\serieslogo{logo_ttl}%please provide filename (without suffix)
\volumeinfo%(easychair interface)
  {M. Antonia {Huertas}, Jo\~ao {Marcos}, Mar\'ia {Manzano}, Sophie {Pinchinat}, \\
  Fran\c{c}ois {Schwarzentruber}}% editors
  {5}% number of editors: 1, 2, ....
  {4th International Conference on Tools for Teaching Logic}% event
  {1}% volume
  {1}% issue
  {191}% starting page number
\EventShortName{TTL2015}
\begin{document}

\maketitle

\begin{abstract}
I have developed a pedagogy and textbook for teaching logic centered on what I call “logical worldviews”.  A logical worldview examines the close connection between philosophical commitments and the logical principles and method for a particular historical logical system.  The class examines multiple historical logical worldviews to show how philosophical positions and logical systems have been closely intertwined, and how changes in philosophical positions have over time yielded corresponding changes in the development of logic.  Such an approach has great benefits for teaching logic to undergraduates.
\end{abstract}

\section{A \textquotedblleft New\textquotedblright\ Approach to Teaching
Logic to Undergraduates: Logical Worldviews}

Since 2006, I have developed an approach to teaching logic to undergraduate
philosophy majors that bases logic instruction in the philosophical
background underpinning historic logical systems. Starting in 2013, I
expanded this pedagogical approach to a general education introductory logic
class, which is open to students of any major and at any stage of their
college career. This curricular development has culminated in a textbook
designed for an introductory logic course that can be used by teachers with
minimal or no experience in the history of logic.

In the textbook I present a re-contextualization of classical formal logic
as just one of many approaches to logic, each having its own philosophical
commitments and unique value. I suggest that there is more to logic than is
construed by only considering its contemporary formal aspects; after all,
formalism in mathematics and logic is a relatively recent innovation,
adopted in part because of other, non-mathematical considerations and
commitments.

I hold that the history of logical thought includes a series of logics, or
logical systems, each of which has non-formal philosophical commitments. I
further hold that these commitments are not separate and independent add-ons
to the logical system, but that the logical method and fundamental
principles of the system are shaped and determined by these commitments.
Though a historical logic can be considered formally and in the abstract,
such a perspective yields an inaccurate account of what the logical system
was intended to do by its author, of the true power of its methods of
reasoning, and of the precise meaning of the logical principles employed.
When the logical system is considered in conjunction with the intertwined
philosophical commitments, I characterize this broader conceptualization as
a \textquotedblleft logical worldview.\textquotedblright 

\section{Using Philosophy to Better Understand a Logical System}

It will be useful at this point to consider a particular example, and
Aristotle is perhaps the best place to start: after all, as the first to
fully develop a logical worldview, many later logics were (as I present in
the textbook) a response to rejecting one or more Aristotelian philosophical
commitments. What follows is a very brief (and incomplete) characterization
of some of the philosophical commitments in Aristotle's logic.\footnote{%
Here and throughout this paper, I present the philosophical views in an
introductory and casual manner, to match the form in which they are
presented in the textbook. Though some of these Aristotelian interpretations
are contested by scholars of Aristotle, a close and well-defended history of
Aristotle's logical and philosophical positions would not be to the point.
What matters for the sake of this pedagogical approach is that
philosophically and historically inexperienced students see that Aristotle's
logical system is amenable to being better understood by framing it in terms
of philosophical positions that Aristotle (plausibly) held.} 

\begin{itemize}
\item Theory of truth: Correspondence. Propositions are true if and only if
they correspond to how the world actually is; false otherwise. Somewhat
anachronistically, we might refer to propositions as being true if and only
if they correspond to facts about the world. This view also includes
bivalence -- propositions are either true or false.\footnote{%
See Aristotle \emph{Categories} chapter 14.} 

\item Access to truth: Direct Realism. Human beings can directly access the
world, and the properties and entities they perceive are (generally, in
principle) the same as the actual properties and entities that exist in the
world.\footnote{%
See Aristotle \emph{de Anima}, particularly Book II chapter 12 through Book
III chapter 4.} 

\item Theory of mind: Simple Apprehension and Judgment. The human mind has a
capacity to discern the true properties of things. Further, we are capable
of abstracting from particular properties to universal forms shared by more
than one individual entity, and recombining such forms in the imagination
and in memory. Finally, we can combine previously unseen combinations of
forms and judge whether or not such propositions correspond to reality.%
\footnote{%
See Aristotle \emph{On Interpretation }chapter 1 and Posterior Analytics
Book II chapter 19.} 

\item Metaphysics: Substances, Forms and a Hierarchical Ontology. In
Aristotle's account, the world is constituted by individual entities, or
primary substances, each of which has a number of properties, or forms.
Substances can be grouped together into secondary substances, or species,
for which some forms are essential, others not. Species can be grouped
together into a higher-order grouping, or genus, which have their own
essential forms. Collectively the interconnection of species and genus
through essential forms constitutes a hierarchical ontology, which is
rationally discernible.\footnote{%
See Aristotle \emph{Categories}; also \emph{Metaphysics.}} 
\end{itemize}

In my textbook I explain how these (and other) features of Aristotle's
philosophical worldview determine the specifics of his formal logic. I'll
sketch this dependence briefly here.

\bigskip 

First, individual substances and their forms are the fundamental perceivable
elements of reality. Collectively, these are referred to as
\textquotedblleft terms\textquotedblright , which are the atoms of
Aristotle's logical system.

Second, propositions are assertions about the actual facts of the world. Any
fact will fundamentally involve two components: an identified subject and
one of its forms. Since propositions correspond to these fundamental facts,
every proposition will have only two terms: a subject term (what the
proposition is about), and a predicate term (what is being said of that
subject).

Third, through experience, induction, abstraction, and natural talent,
humans are capable of discerning some of the essential qualities of
substances. This allows us to generalize beyond particular substances to be
able to quantify the scope of the assertion (all, some, one) and whether the
subject and predicate are linked (affirmation) or separated (denial).

Fourth, through our capacity to discern the ontological structure of
reality, we are able (in principle) to choose any two propositions and
connect them by means of something in the hierarchy that links them (a
middle term). The fundamental form of reasoning from one proposition (with
two terms) and another (sharing one of those terms) is through a bridge
proposition, sharing a \textquotedblleft middle\textquotedblright\ term with
our premise and conclusion. This is the minimum size of an argument
productive of new knowledge, and thus the syllogistic form must have three
propositions sharing exactly three distinct terms.

This is a rough and quick sketch of how Aristotle's philosophical
commitments determine the form of his logical method, which culminates in
the syllogism. This perspective on Aristotle's logical worldview also helps
to explain many other curiosities about Aristotle's logic that have
perplexed many trying to understand his logic. I can't get into them in
detail here, but here are a few positions in Aristotle's logic that are
explained by considering his philosophical commitments.

\begin{itemize}
\item Reduction to the first figure. Rather than be satisfied with
identifying valid syllogistic moods, he takes great pains to reduce them to
four figure one moods. He does so because he believes that these four moods
are the fundamental natural and self-evident reasoning patterns for human
beings.\footnote{%
See Aristotle \emph{Prior Analytics} Book II chapters 1-15.} 

\item The missing fourth figure and two moods. In inventorying valid moods,
Aristotle neglects to identify two valid moods in figure one and figure four
altogether. These moods do not fit the \textquotedblleft
natural\textquotedblright\ reasoning patterns determined by Aristotle's view
and are thus neglected as inappropriate forms of human reason.\footnote{%
For a good discussion of this issue see Lynn E. Rose, \textquotedblleft
Aristotle's Syllogistic and the Fourth Figure,\textquotedblright\ \emph{Mind 
}vol. 74, no. 295 (July 1965): 382-389.} 

\item Existential import. Many have found that Aristotle's inference that
\textquotedblleft All S are P\textquotedblright\ implies \textquotedblleft
Some S is P\textquotedblright\ to be logically invalid. However, given his
philosophical commitments, no universal claim could ever be made without
first establishing the truth of the particular, so the truth of the
universal will always imply the particular.
\end{itemize}

Though this exposition of this logical worldview is seriously abridged, it
should demonstrate that the form, method and logical principles of
Aristotle's logic are heavily dependent on Aristotle's philosophical
commitments. Those who wish to teach Aristotle's logic to students can use
the logical worldview in order to have Aristotle's logic make much more
sense to students, and show how his logic is closely connected to many of
Aristotle's other philosophical positions.

\section{Using Philosophical Commitments to Explain Transitions from One
Logical System to Another}

I'll briefly discuss how later logical systems can be understood as having
rejected or modified some part or other of Aristotle's logical worldview.
I'm not suggesting that these later logics were the result of simply
abandoning some philosophical commitment or other of Aristotle's; rather, by
rejecting one or more commitments, new opportunities for logical method are
opened up, which in combination with some innovations led to new ways to do
and conceive logic. I'll briefly consider here two logical worldviews
post-Aristotle.

\bigskip 

\emph{Bacon's rejection of Aristotle's theory of mind: a motivation for
natural history and modern natural sciences.} In Bacon's \emph{New Organon},
he takes Aristotelian philosophers to task for relying too heavily on our
natural capacity to understand the forms present in reality and argues that
nature is far too subtle for our unaided sensory and cognitive capacities.
He argues that we must develop new instruments to enhance our senses
(microscopes, thermometers, etc.) and to enhance our cognitive discernment
of what we are seeing by means of experiments.\footnote{%
He makes this case over the course of Book I.} Since our capacity to
naturally discern the forms in reality is flawed, we must develop a new
logical technique, which is basically induction from a very large number of
particular instances towards generalizations of increasing universality and
confidence.\footnote{%
This technique is laid out in Book II.} Therefore, the syllogism, which
relies on gathering information about forms and then deducing new truths is
for Bacon the wrong logical method. Bacon does not reject Aristotle's other
philosophical commitments mentioned above, and to varying degrees and with
some modification, the natural sciences still endorse a correspondence
theory of truth, a capacity to discern forms of substances, and a
hierarchical ontology.

\bigskip 

\emph{Boole's rejection of Aristotle's Direct Realism and metaphysics: an
effort to mathematize logic.} To a great extent Newton's mathematical
approach to physics follows similar rejections as Boole, but for the sake of
what will soon become clear, Boole's contribution to logic is easier to
present to students. Starting with the Early Modern philosophers in the 17th
century and culminating in the modern approach to science in the 19th, the
rejection of Direct Realism became fairly common. Philosophers were no
longer inclined to believe that we had direct access to the substances and
forms through our senses exactly as they are in reality. Because we lacked
the capacity to see the world as it actually was, Boole changed the focus of
logic from terms to classes in his seminal work The Laws of Thought. Instead
of the logical atom being a term, which required discerning the true forms
of individuals, he made the fundamental element the class, a collection of
individuals gathered in some (perhaps unknowable, perhaps arbitrary) way.
One could identify a class without having to know exactly what was required
for membership. His particular brilliance was to show that the logical
method of Aristotle could be encompassed within an algebra for classes. By
shifting from term to class, Boole needed a new logical method, which he
found in algebra. Though Boole mathematized logic, one should note that he
did not abandon Aristotle's philosophical commitments completely: he still
largely conceived of reality as being metaphysically ordered in a hierarchy,
and also that truth was established by correspondence. Further, he did
believe that he had discovered the \textquotedblleft laws of
thought\textquotedblright : the actual laws by which human minds reasoned,
maintaining a close link between logic and theory of mind.\footnote{%
See Boole \emph{The Laws of Thought} chapter 1. As an interesting
philosophical divergence into the history of computing and artificial
intelligence, I often have students learn and discuss the development of
logic machines in the 1870s following Boole's logical system. These machines
prompted a brief but lively debate amongst American logicians such as C. S.
Peirce as to whether machines that followed the laws of thought actually
were thinking. (see C. S. Peirce, \textquotedblleft Logical
Machines,\textquotedblright\ \emph{The American Journal of Psychology}, vol.
1, no. 1 (Nov. 1887): 165-170) The view seemed generally to be in the
affirmative. Eventually, however, those doing algebra of logic dropped their
commitment to linking logic to laws of thinking. (see, for example, chapter
1 in Louis Couturat, \emph{The Algebra of Logic} (Chicago: Open Court
Publishing, 1914))} 

\bigskip 

I hope that these two brief sketches outline how disparate logical systems
in history may appear to be radically different, but that their differences
are in part due to rejecting or modifying some basic philosophical (and
non-logical) positions. These transitions can show students that historical
logics are not useless or obsolete; instead, they are evolutionary
predecessors to later logics (and ultimately contemporary classical logic).%
\footnote{%
By the way, in my textbook I use logical worldviews to show how contemporary
propositional logic results from abandoning a commitment to a correspondence
theory of truth and the resulting modification of the algebraic method of
Boole's logic. This philosophical position is partly inspired by formalism
-- often identified with Giuseppe Peano and many of his Italian colleagues.
See Paolo Mancosu, Richard Zach, and Calixto Badesa, \textquotedblleft The
Development of Mathematical Logic from Russell to Tarski,
1900-1935\textquotedblright\ in \emph{The Development of Modern Logic},
edited by Leila Haaparanta, p. 318-470. Oxford: Oxford University Press,
2009.} These historical logics can still be valuable and useful to the
contemporary student, should one embrace the philosophical commitments
underlying the logical worldview. 

Over the years, I have taught and analyzed a number of different historical
logical systems and perspectives by means of a close examination of their
philosophical commitments. In addition to those systems mentioned above, I
have found that at least the following logical approaches can be analyzed in
this way: Medieval faith and reason, propositional and predicate logic,
mathematical physics (a la Newton), Ancient/Medieval problem of universals,
nominalism, Leibniz's logical innovations, logical fatalism, logical
atomism, logicism, and the theory of relations. I anticipate many more, if
not most, historical logics and logical issues would also be amenable to
this pedagogical approach.

\section{Pedagogical Benefits of Teaching Logical Worldviews}

In teaching introductory formal logic to undergraduates one may have
encountered the following scenario: those students who have an affinity to
formal reasoning (e.g. math, computer science, and science students) take
rather well to logic, but other students (e.g. humanities and social science
students) struggle. This should not be too surprising, given the
mathematical basis for contemporary formal logic. But what to do with those
students left behind? And though the power and flexibility of formal logic
is well demonstrated, especially relative to earlier approaches to logic
(such as Aristotelian logic), students that are not mathematically inclined
find formal logic both confusing and limiting. Such students may be capable
of reasoning well, and may have a great deal of experience and success in
reasoning, but often find that their previously successful reasoning
practices are not well captured by the methods found in classical formal
logic. One solution to this mis-fit is addressed by courses in informal
logic (sometimes called \textquotedblleft critical
thinking\textquotedblright ). And though such an approach to reasoning more
closely fits the more natural, intuitive, and practical forms of reasoning
of interest to many college students, it typically achieves this result by
minimizing or avoiding formal reasoning altogether.

The contemporary setting of formal logic in introductory undergraduate
education is a difficult, if not paradoxical one. Students not already open
to formal reasoning often find formal logic not very \textquotedblleft
useful\textquotedblright , and many struggle to see the value or purpose of
learning it. This is further compounded when formal logic is taught by
philosophy departments, a discipline that otherwise frequently engages in
applied reasoning in a variety of fields of contemporary relevance (science,
religion, art, ethics, etc.), and in \textquotedblleft
big\textquotedblright\ questions that are quite meaningful to students in
their lives (\textquotedblleft What is the meaning of
life?\textquotedblright\ \textquotedblleft What is the best life one can
live?\textquotedblright\ \textquotedblleft What is the right thing to
do?\textquotedblright ). It is rare for a logic class to turn from the study
of formal reasoning methods to connecting those to issues of current
interest (\textquotedblleft informal\textquotedblright\ or philosophical
reasoning), and so students who do ultimately learn formal logic have great
difficulty applying their new skills to questions of particular interest to
them, unless they continue on with more advanced philosophy courses.

By understanding logical systems and methods in the context of their
philosophical commitments, students are exposed to the close connection
between logic and philosophy. They learn why logic has the method it does,
which helps them engage and understand the material more deeply than simply
challenging them to master the mathematical methods of contemporary formal
logic. Further, students learning logical worldviews are more capable of
knitting together formal logic and \textquotedblleft
informal\textquotedblright\ or philosophical reasoning, as these approaches
are explicitly linked. Finally, logical worldviews demonstrate how logic is
closely connected to philosophy, and it is rather easy to introduce
philosophical questions of interest to students in the midst of learning a
logical system.\footnote{%
For example, students can explore logical fatalism (and their more
philosophically interesting versions concerning free will and determinism or
divine foreknowledge) by exploring the consequences of Aristotle's
commitment to a correspondence theory of truth and bivalence. Also, students
can better understand Anselm's ontological argument and Guanilo's reply by
seeing how Anselm's position is grounded in Aristotle's own reductio proofs,
and Guanilo's on Aristotle's theory of mind and direct realism.} 

\bigskip 

\section{Benefits of Teaching Logical Worldviews to History of Philosophy}

Most logic textbooks, and by implication, most undergraduate logic classes
teach a system of formal logic developed in the 20th century (except, in
some cases, a cursory examination of syllogistic reasoning). This approach
to logic has some disturbing implications for the 2,500 years of
philosophical reasoning that preceded the 20th century. The considerate
student (as well as many professors) may draw any or all of the following
conclusions about pre-20th century logical reasoning:

\begin{enumerate}
\item Logic from these time periods was either wrong, incomplete, or
unacceptably limited in scope and methods. After all, if such logics were
suitable, why develop contemporary formal logic that encompasses and
surpasses such logics?

\item By extension from this first point, philosophical reasoning based on
earlier logics could similarly be considered to be flawed, obsolete, or
otherwise useless. Some schools of philosophical thought have explicitly
embraced this perspective (e.g. the Logical Positivists); today, this
perspective is often implicit in attitudes about historical philosophy from
a wide variety of fields (every so often one will make such an assertion
explicitly, as Stephen Hawking recently did).\footnote{%
See Stephen Hawking and Leonard Mlodinow, \emph{The Grand Design} (New York:
Bantam Books, 2010) chapter 1: \textquotedblleft Traditionally these are
questions for philosophy, but philosophy is dead. Philosophy has not kept up
with modern developments in science, particularly physics. Scientists have
become the bearers of the torch of discovery in our quest for
knowledge.\textquotedblright } 

\item Many contemporary philosophers (students and professors) educated in
contemporary formal logic enough to use it to parse examples of natural
reasoning will use it to parse examples of philosophical reasoning that
pre-dates mathematical logic. This is anachronistic to say the least, and
results in a misleading conception of the actual reasoning made by these
historic figures. Further, it is likely to construe historical philosophy in
a bad light, since it makes it appear that these historical figures engaged
regularly in invalid reasoning for most of the intellectual history of
Western thought.
\end{enumerate}

For some, these considerations are not of much concern: many do, in fact,
believe that history is rife with examples of bad reasoning, just as up to
the scientific revolution, history was plagued with examples of bad science.
This perspective is reinforced by the (mostly legitimate) observation that
some of history's most highly regarded philosophical reasoners seem to have
made grave errors in their knowledge claims (Aristotle is a significant
target here).

In math and science an ahistorical perspective may be appropriate: after
all, most science and math textbooks teach the completed current state of
knowledge, not the long and difficult (and in many cases wrong) path that
led us here. Outside of these fields, however, the value of accurately
understanding past perspectives is more apparent: many areas of study in the
humanities and social sciences to this day engage well-reasoned theories
dating back to Ancient Greece.

Setting concerns about history of philosophy aside, a large number of people
presently endorse Aristotle's philosophical commitments enumerated above and
teachers would likely find that most of those people would have little
problem with syllogistic reasoning. In contrast, I expect that most teachers
of logic find that many students find much of contemporary formal logic
counter-intuitive and alien to how they reason. In contrast, logical
worldviews can help the logic teacher show that in certain philosophical
contexts, particular historic logics may be the best approach; but that
should those philosophical commitments be rejected or disproved in other
reasoning contexts, different logics are called for. A student's intuitive
reasoning can be embraced rather than replaced, which would certainly seem
to be a desirable outcome for a logic class.

\section{Benefits for the Teacher of Logic}

I'd like to finish up by discussing a few benefits for the logic teacher who
adopts this approach. I've designed the textbook to contain a number of
different modules, each examining a different logical worldview. As I
complete more modules, the textbook will fill up with a fair number of
unique instances of logical systems spanning the last 2,500 years. Each
module contains several key components: (1) an inventory of the
philosophical commitments of the logical system; (2) a trace of the
connection between the philosophical commitments and the development and
form of the logical method and principles employed; (3) an introduction and
discussion of some interesting philosophical issues that arise under this
worldview.

Each module is designed to be teachable by those with only a minimal
expertise in the history of philosophy.\footnote{%
As mentioned earlier, this involves presenting a simplistic interpretation
which historians of philosophy might find unsupported by sufficient research
and scholarship. However, most textbooks are likely subject to the same
complaint, and thus this publication genre permits, if not requires, such
oversimplification.} The characterization of a logical worldview does depend
on a historical interpretation of the philosophy and logic under study, but
a background in historical scholarship is not required in order to
understand or teach the main components of each module. Instead, anyone with
a good background in philosophy should be able to comprehend and easily
teach the first and third components (inventory of philosophical
commitments, and linkage to philosophical issues). The second component of
each module (tracing the connection between the philosophy and the logic) is
more subtle and complex, but should be within reach for anyone experienced
in teaching logic.

Modules are designed to be able to be taught independently, although many
are closely related to others. This allows the teacher to pick and choose a
path through the history of logic, and thereby emphasize the impact of
different philosophical views on the development of logic and/or our
understanding of certain logical concepts and methods. For example, in a
logic class that culminates in teaching propositional and predicate logic,
one might start with Aristotle and then pick modules that progressively show
how logic became more mathematical and abstract (by for instance, selecting
logical worldviews of Leibniz, Boole and Russell as transitional modules).
On the other hand, a course interested in the logic of the scientific method
might feature the logical worldviews of Aristotle, Medieval faith and
reason, Bacon, and Newton.

In conclusion, I have found this approach to be of great benefit to students
and to myself by expanding teaching of logic beyond the boundaries of formal
methods and systems. Students seem to derive a deeper understanding of logic
by seeing how it connects with other philosophical positions. Further, this
approach helps make logic more comprehensible to the broad pool of students
who are either disinclined or not well suited towards the strictly formal,
mathematical approaches to teaching logic so commonly taught today.

\end{document}